\documentclass[a4paper,11pt]{article}
\usepackage[utf8]{inputenc}

\usepackage{geometry}
\usepackage{tikz}
 \usepackage{pgfplots}
 \pgfplotsset{compat=1.3}
\usepgfplotslibrary{groupplots}

\usepackage[english]{babel} 
\usepackage{amssymb}
\usepackage{amsmath}
\usepackage{amsthm}
\usepackage{graphicx}
\usepackage{txfonts}
\usepackage{color}
\usepackage{subfigure}
\usepackage{multirow}
\usepackage{hyperref}

\newtheorem{assumption}{Assumption}
\graphicspath{ {Images/} }
\renewcommand*{\vec}[1]{\boldsymbol{#1}}

\newdimen\nodeDist
\tikzset{
    position/.style args={#1:#2 from #3}{
        at=(#3.#1), anchor=#1+180, shift=(#1:#2)
    }
}

\title{\bf {A computational study of preconditioning techniques for the stochastic diffusion equation with lognormal coefficient}}
\author{ Eugenio Aulisa \footnote{Department of Mathematics and Statistics, Texas Tech University, Lubbock TX 79409, USA.}%
  \and  Giacomo Capodaglio \footnote{Department of Scientific Computing, Florida State University, Tallahassee FL 32306, USA.}%
  \and Guoyi Ke \footnote{Department of Mathematics and Physical Sciences, Louisiana State University at Alexandria, Alexandria LA 71302, USA.}%
  }
\date{}

\begin{document}

\maketitle

\begin{abstract}
We present a computational study of several preconditioning techniques for the GMRES algorithm applied to the stochastic diffusion equation with a lognormal coefficient discretized with the stochastic Galerkin method. 
The clear block structure of the system matrix arising from this type of discretization motivates the analysis of preconditioners designed according to a field-splitting strategy of the stochastic variables. 
This approach is inspired by a similar procedure used within the framework of physics based preconditioners for deterministic problems, and its application to stochastic PDEs represents the main novelty of this work.
Our numerical investigation highlights the superior properties of the field-split type preconditioners over other existing strategies in terms of computational time and stochastic parameter dependence.
\end{abstract}


\section{Introduction}
In the last decade, stochastic partial differential equations (SPDEs) have attracted great attention from the scientific community, 
due to their ability to take into account uncertainties entering the problem through the input data. 
These sources of uncertainty may arise for instance from boundary and initial conditions, coefficients, forcing terms, 
or intrinsic randomness of the processes as in the case of heterogeneous media 
\cite{le2003multigrid,ghanem1998probabilistic, matthies1999finite, cliffe2011multilevel, xiu2005high, nobile2008sparse}.
The solutions of SPDEs allows to characterize the mean, 
variance and in general the probability density function of quantities of interest in the post-processing phase \cite{capodaglio2019piecewise, capodaglio2018approximation, powell2009block}.
Among a wide variety of numerical methods for solving SPDEs \cite{gunzburger2014stochastic},
a popular approach is the stochastic Galerkin method (SGM), 
where a Galerkin projection is employed to approximate the infinite dimensional stochastic space with a finite dimensional one, spanned by appropriate basis functions \cite{ghanem2003stochastic}.

In this work, we focus on the efficient numerical solution of the stochastic diffusion equation with a lognormal coefficient, 
arising for instance in the framework of groundwater flows, where the permeability coefficient is often considered to be lognormal \cite{freeze1975stochastic, lu2002stochastic}. 
Using the SGM, the resulting stiffness matrix is block dense, due to the nonlinearity of the coefficient \cite{ullmann2015solving}, and ill-conditioned with respect to the mesh size and to the stochastic parameters such as the standard deviation of the input lognormal field \cite{ullmann2012efficient, ernst2010stochastic, powell2010preconditioning, ullmann2010kronecker}. The aforementioned properties of the matrix require the design of ad hoc preconditioned solvers for the efficient solution of the SGM system.
When the stochastic diffusion problem is well-posed, regardless of the type of coefficient, the system matrix is symmetric and positive definite, hence a huge variety of preconditioned conjugate gradient (PCG) solvers has been proposed in the literature \cite{desai2018scalable, benner2015low, subber2014schwarz, brezina2014smoothed, ullmann2012efficient, rosseel2010iterative, powell2009block, rosseel2008algebraic}.
As pointed out in \cite{ullmann2015solving} however, the density of the matrix given by a lognormal coefficient makes the use of PCG methods problematic, given that for every PCG iteration a matvec operation has to be performed.
On the other hand, in two recent studies very relevant to our framework \cite{ullmann2012efficient, ullmann2015solving}, the preconditioned GMRES algorithm has been shown to perform better than PCG in terms of both solution time and dependence on the stochastic parameters for a stochastic diffusion equation with lognormal coefficient. In the above mentioned studies, the diffusion problem is reformulated as a convection-diffusion problem with the result that the nonlinear coefficient is transformed into a linear one, and the resulting matrix gains sparsity, while losing symmetry. The preconditioned GMRES from \cite{ullmann2012efficient} and \cite{ullmann2015solving} showed independence on spatial discretization and most stochastic parameters, however a mild but relevant influence on the standard deviation of the input lognormal field was observed (with the number of iterations going from 6 to 27 in the worst case scenario in \cite{ullmann2012efficient}). Hence, the important work done in the two studies described above motivated our choice of studying the performances of a preconditioned GMRES algorithm rather than a PCG as most studies in the literature have done.
We also mention \cite{jin2007parallel} for a work that has employed a flexible GMRES algorithm as a solver.

After choosing the solver, the next crucial step is the choice of appropriate preconditioners. 
Several different strategies have been carried out for this task, with the most popular possibly being the so called mean-based preconditioner \cite{ghanem1996numerical, powell2009block}. Variants of the mean-based preconditioner have also been designed \cite{ullmann2010kronecker}.
Moreover, domain decomposition type methods have also been used as preconditioning techniques \cite{subber2014schwarz, desai2018scalable, jin2007parallel}, as well as low-rank approaches \cite{benner2015low} and other approaches such as multigrid \cite{rosseel2010iterative, rosseel2008algebraic}.
In particular, the extensive study on iterative solvers for  stochastic PDEs carried out in \cite{rosseel2010iterative} reported that for the case of a lognormal field, only a CG solver preconditioned with multigrid with a block Gauss-Seidel smoother showed robust convergence, as the other approaches considered suffered a dependence on the stochastic parameters, such as the variance of the input field.
The study in \cite{rosseel2010iterative} also concluded that for large problems, multigrid type methods should be preferred. Hence, in light of \cite{rosseel2010iterative}, in the present study we decided to consider preconditioners of multigrid type, specifically geometric multigrid.
For completeness, we also mention works that employed geometric or algebraic multigrid as a solver, although not in the context of a lognormal diffusion coefficient \cite{le2003multigrid, elman2007solving}.

The major novelty of this work is the introduction of block preconditioners for the multigrid smoother with a structure arising from using the stochastic modes as a field-splitting (FS) strategy. This approach is motivated by promising results on FS preconditioners for deterministic PDEs obtained by the authors in a series of papers  \cite{ke2017block,ke2018augmented,ke2018new, calandrini2019field}.
When used on deterministic PDEs, the FS strategy yields a block structure associated with the physical variables on the physical domain. 
On the other hand, here the stochastic modes are used in an analogous way as physical variables for the splitting strategy.
The FS approach can also be applied directly as a preconditioner for the GMRES algorithm although, as it will be shown, the best computational performances are obtained if FS is used within the framework of geometric multigrid, i.e. on the smoother. To the best of our knowledge, this is the first work to perform a computational analysis of the performances of GMRES preconditioned with geometric multigrid in the framework of SPDEs. Because from now on we will be dealing only with geometric multigrid, MG will implicitly refer to its geometric variant.

The paper is organized as follows: 
Section \ref{numMod} presents the mathematical formulation of the stochastic diffusion problem with lognormal coefficient and its weak form, discretized with the stochastic Galerkin method. 
In Section \ref{preconds}, five preconditoned GMRES solvers are introduced, including those using FS, and their implementation is discussed.
In Section \ref{numres}, the performance of the solvers introduced in the previous section is assessed through a thorough numerical study, focusing on computational time,  number of iterations, 
and their dependence on the spatial parameter and stochastic parameters.
Finally, conclusions are drawn in the last section.

\section{Numerical Modeling}\label{numMod}
Let $(\Omega, \mathcal{F}, \mathbb{P})$ be a probability space, where $\Omega$ is the set of outcomes,
$\mathcal{F} \subset 2^\Omega$ is the $\sigma$-algebra of events, and $\mathbb{P}: \mathcal{F} \rightarrow [0,1]$ is a probability measure.
We define the spatial domain $D \subset \mathbb{R}^d, d =1, 2, 3$ and denote its boundary with $\partial D$.
Then, a general stochastic diffusion problem reads: find  $u:\overline{D} \times \Omega \rightarrow \mathbb{R}$ 
such that the following equations hold $\mathbb{P}$-almost everywhere in $\Omega$

\begin{align} 
  -\nabla \cdot \big(a(\vec{x},\omega) \nabla u(\vec{x},\omega)  \big) &= f(\vec{x}), & \mbox{ in } D \times \Omega, \label{Strong_PS}\\
  u(\vec{x},\omega)         				    &= 0, & \mbox{ on } \partial {D} \times \Omega. 
\end{align}

Note that, for simplicity,  we only consider the stochastic contribution coming from the coefficient function and not from the forcing term.
To guarantee that the problem in \eqref{Strong_PS} is well-posed, we make the following assumptions, similar to \cite{gunzburger2014stochastic, babuvska2007stochastic, nobile2008sparse, nobile2008anisotropic, capodaglio2018approximation}, 
\begin{assumption} \label{ass1}
The coefficient function $a(\vec{x},\omega)$ in system \eqref{Strong_PS} has the following properties: \\
1. There exists a positive constant $ a_{min} < \infty$ such that $a_{min} \leq a(\vec{x},\omega) $ almost
surely on $\Omega$, for all $\vec{x} \in \overline{D}$. \\
2. $a(\vec{x},\omega)= a(\vec{x},\vec{y}(\omega))$ in $\overline{D} \times \Omega$, 
where $\vec{y}(\omega) = (y_1(\omega), y_2(\omega), ..., y_N(\omega))$ is a vector of real-valued uncorrelated random variables. \\
3. $a(\vec{x},\vec{y}(\omega))$ is measurable with respect to $\vec{y}$. 
\end{assumption}

For any $n = 1, \dots, N, $ let $\Gamma_n :={y_n(\Omega)} \in \mathbb{R}$, and denote the probability density function (PDF) of $y_n$ by
$ \rho_n(y_n ): \Gamma_n \rightarrow \mathbb{R}$. Then the joint image of ${\{y_n(\omega)\}_{n=1}^N}$ can be defined as $\Gamma := \bigcup_{n=1}^{N}\Gamma_n$ with 
the joint PDF $ \rho(\vec{y} ) : \Gamma \rightarrow \mathbb{R}$.   

\subsection{Lognormal random field}
Define the random field $ \gamma(\vec{x},\vec{y}(\omega))$ as,
\begin{equation}
 \gamma(\vec{x},\vec{y}(\omega)) := \textrm{log}(a(\vec{x},\vec{y}(\omega)) - a_{min}),
\end{equation}
where $a(\vec{x},\vec{y}(\omega))$ is the random field in Eq. \eqref{Strong_PS} satisfying Assumption \ref{ass1}.
The truncated Karhunen-Lo\`eve (KL) expansion \cite{li2008fourier, schevenels2004application, frauenfelder2005finite} is used to approximate $\gamma$, yielding
\begin{equation}\label{aKL}
\gamma(\vec{x},\vec{y}(\omega)) \approx \gamma_{KL}(\vec{x},\vec{y}(\omega)) := \mu_{\gamma} +\sum_{n=1}^{N} \sqrt{\lambda_n} b_n(\vec{x})y_n(\vec{y}(\omega))
\end{equation}
where $\mu_{\gamma}$ is the mean of $\gamma(\vec{x},\vec{y}(\omega))$, ($\lambda_n$, $b_n(\vec{x}$)) 
is the $n^{th}$ eigenpair of the covariance function of $\gamma(\vec{x},\vec{y}(\omega))$ whose eigenvalues $\lambda_n$ are positive,
and listed in non-increasing order. 
It follows that the stochastic coefficient $a(\vec{x},\vec{y}(\omega))$ can be rewritten as
\begin{equation}
 a(\vec{x},\vec{y}(\omega)) = a_{min} + \textrm{exp}(\gamma(\vec{x},\vec{y}(\omega))).
\end{equation}
and approximated by 
\begin{equation} \label{KL_gamma}
 a(\vec{x},\vec{y}(\omega)) \approx a_{KL}(\vec{x},\vec{y}(\omega)) = a_{min} + \exp \bigg(\mu_{\gamma} +\sum_{n=1}^{N} \sqrt{\lambda_n} b_n(\vec{x})y_n(\omega)\bigg).
\end{equation}
Throughout the paper we assume that $\vec{y}$ is Gaussian, hence the random variables ${\{y_n(\omega)\}}_{n=1}^N$ are standard independent and identically distributed. It follows that the random field $\gamma(\vec{x},\vec{y}(\omega))$ is Gaussian and so $a(\vec{x},\vec{y}(\omega))$ is lognormal.
The eigenpair ($\lambda_n$, $b_n(\vec{x}$)) is obtained with the solution of the following generalized eigenvalue problem
\begin{equation}\label {gen_eigen}
\int_{D} C_{\gamma} (\vec{x},\vec{\hat{\vec{x}}}) b_n({\vec{x}}) d \vec{x} = \lambda_n b_n({\hat{\vec{x}}}),
\end{equation}
where $C_{\gamma} (\vec{x}, \hat{\vec{x}})$ is the covariance function of the field $\gamma(\vec{x},\vec{y}(\omega))$. 
Here the covariance function of $\gamma(\vec{x},\vec{y}(\omega))$ is assumed to be
\begin{equation} \label{standard_deviation}
 C_{\gamma} (\vec{x},\hat{\vec{x}}) = \sigma_{\gamma}^2 \exp \bigg[- \frac{1}{L_c} \Big(\sum_{i=1}^d |x_i - \hat{x}_i|\Big)\bigg],
\end{equation}
where $\sigma_{\gamma}$ denotes the standard deviation of $\gamma(\vec{x},\vec{y}(\omega))$, 
$d$ is the dimension of the spatial variable and $L_c>0$ is a correlation length satisfying $L_c \leq \textrm{diam}(D)$.  
Details on how to solve Eq. \eqref{gen_eigen} using a Galerkin procedure can be found in \cite{capodaglio2018approximation}.
For ease of notation, from now on we will drop the $\omega$ from all $\vec{y}(\omega)$ and $y_n(\omega)$. 
\subsection{The Stochastic Galerkin Method}
The weak formulation of the stochastic diffusion problem consists of finding  
$u(\vec{x},\vec{y})\in \mathbb{W} = H_0^1(D) \bigotimes L^2(\Gamma)$ satisfying 
\begin{equation} \label{weak_form}
\int_{\Gamma} \int_{D} a(\vec{x},\vec{y}) \nabla u(\vec{x},\vec{y}) \nabla v(\vec{x},\vec{y}) d\vec{x} d\vec{y}
= \int_{\Gamma} \int_{D} f(\vec{x},\vec{y}) v(\vec{x},\vec{y}) \rho(\vec{y}) d\vec{x} d\vec{y},
\end{equation}
for all $v(\vec{x},\vec{y}) \in \mathbb{W}$.

The above equation can be discretized with the SGM introducing a finite dimensional space to approximate the infinite dimensional space $L^2(\Gamma)$. 
Let $p \in \mathbb{N}$ denote the polynomial order of the associated approximation 
and let $\mathcal{P}_{\mathcal{J}(p)}(\Gamma) \subseteq L^2(\Gamma)$ be a multivariate polynomial space over $\Gamma$ associated with 
the index set $\mathcal{J}(p)$, defined by 
\begin{equation} \label{multivariate_polynomial}
\mathcal{P}_{\mathcal{J}(p)}(\Gamma) = \mathrm{span} \bigg\{ \prod_{n=1}^N y_n^{p_n} \, \Big| \, \vec{p} \in \mathcal{J}(p),y_n \in \Gamma_n \bigg\}{}, 
\end{equation}
with $\vec{p} = (p_1, p_2, \dots, p_N)$. In this work the index set is chosen as
\begin{equation}
 \mathcal{J}(p) = \bigg\{ \vec{p} \in \mathbb{N}^N \,\Big | \, \sum_{n=1}^N p_n \leq p \bigg\}.
\end{equation}
Other choices of $\mathcal{J}(p)$ are possible, see for instance \cite{gunzburger2014stochastic}.
The dimension of $\mathcal{P}_{\mathcal{J}(p)}(\Gamma)$ is finite and is given by $ M_p = (N+p)!/(N! p!)$. 
With the finite element space $V_h(D) \subseteq H_0^1(D) $ and the finite dimensional polynomial space $\mathcal{P}_{\mathcal{J}(p)}(\Gamma) \subseteq L^2(\Gamma)$, 
the approximate solution $u_{hp} \in \mathbb{W}_{hp} =V_h(D) \bigotimes \mathcal{P}_{\mathcal{J}(p)}(\Gamma)$ of the discrete weak formulation satisfies
\begin{equation} \label{discrete_form}
\int_{\Gamma} \int_{D} a(\vec{x},\vec{y}) \nabla u_{hp} (\vec{x},\vec{y}) \nabla v_{hp} (\vec{x},\vec{y}) d\vec{x} d\vec{y}
= \int_{\Gamma} \int_{D} f(\vec{x},\vec{y}) v_{hp} (\vec{x},\vec{y}) \rho(\vec{y}) d\vec{x} d\vec{y},
\end{equation}
for all $v_{hp}(\vec{x},\vec{y}) \in \mathbb{W}_{hp}$.

Fro $V_h(D)$, we choose continuous piecewise-biquadratic polynomials $\{\phi_j(\vec{x})\}_{j=1}^{J_h}$, 
whereas for $\mathcal{P}_{\mathcal{J}(p)}(\Gamma)$, we choose multivariate Hermite polynomials $\{\psi_{\vec{p}}(\vec{y})\}_{\vec{p} \in \mathcal{J}(p)}$.
The multivariate polynomials are obtained in a tensor product fashion from the univariate probabilist Hermite polynomials, which are appropriately scaled so that they form an orthonormal basis with respect to the PDF \cite{capodaglio2018approximation, gunzburger2014stochastic}
\begin{align}\label{univ_pdf}
 \rho_n(y_n) = \dfrac{\exp{(-y_n^2/2)}}{\sqrt{2 \, \pi}}.
\end{align}
This choice is motivated by the assumption of Gaussian distributed input. 
The orthonormal relation means
\begin{equation}
\int_{\Gamma} \psi_{p_i}(y_n) \psi_{p_j}(y_n) \rho_n(y_n) d{y_n} = \delta_{ij},
\end{equation}
where $\delta_{ij}$ is the Kronecker's delta. 
With the two bases introduced, the discrete solution $u_{hp}$ can then be written as
\begin{equation} \label{u_discrete}
 u_{hp}(\vec{x},\vec{y}) = \sum_{\vec{p} \in \mathcal{J}(p)} \sum_{j=1}^{J_h} u_{\vec{p},j}\phi_j(\vec{x}) \psi_{\vec{p}}(\vec{y}) = \sum_{\vec{p} \in \mathcal{J}(p)}  u_{\vec{p} }(\vec{x}) \psi_{\vec{p}}(\vec{y}),
\end{equation}
where we define
\begin{align}\label{u_p}
u_{\vec{p} }(\vec{x}) :=  \sum_{j=1}^{J_h} u_{\vec{p},j}\phi_j(\vec{x}), \quad \vec{u}_{\vec{p} }:= [u_{\vec{p},1}, u_{\vec{p},2}, \dots, u_{\vec{p},J_h}].
\end{align}
The deterministic function $u_{\vec{p} }$ is  the finite element solution corresponding 
to the $\vec{p}^{\textrm{th}}$ stochastic mode and  $\vec{u}_{\vec{p} }$ is the vector of its nodal values. 
Then substituting Eq. \eqref{u_discrete} into Eq. \eqref{discrete_form}, the following linear algebraic system is recovered
\begin{equation} \label{linear_algebraic}
\sum_{ \vec{p'}\in \mathcal{J}(p)}  \bigg( \int_{\Gamma} A(\vec{y})  \psi_{\vec{p}} (\vec{y}) \psi_{\vec{p'}}(\vec{y}) \rho(\vec{y}) d\vec{y} \bigg)
 \vec{u}_{\vec{p'}} = \int_{\Gamma} \vec{f} \psi_{\vec{p}}(\vec{y}) \rho(\vec{y}) d \vec{y},
\end{equation}
where $A_{j,j'}(\vec{y}) =  \int_D a(\vec{x},\vec{y}) \nabla \phi_j(\vec{x}) \cdot \nabla \phi_{j'}(\vec{x}) d\vec{x}$ 
and $\vec{f}_j = \int_D f(\vec{x}) \phi_j(\vec{x}) d\vec{x}$, for $j, j' = 1, \dots, J_h$.

As in \cite{gunzburger2014stochastic, ullmann2012efficient}, we use a generalized polynomial chaos expansion \cite{xiu2002wiener} for the coefficient $a_{KL}(\vec{x}, \vec{y})$ in Eq. \eqref{KL_gamma}: for any $q \in \mathbb{N}$, it is computed as 
\begin{equation} \label{dataProjection} 
 a_{KL}(\vec{x}, \vec{y}) = \sum\limits_{\vec{q} \in \mathcal{J}(q)} a_{\vec{q}}(\vec{x}) \psi_{\vec{q}}(\vec{y}),
\end{equation}
where $ a_{\vec{q}}(\vec{x}) = \int_{\Gamma} a_{KL}(\vec{x}, \vec{y}) \psi_{\vec{q}}(\vec{y}) \rho(\vec{y}) d\vec{y}$. 
Note that the larger $q$ is, the more accurate the projection of the random field will be, however the resulting system matrix will have a denser sparsity pattern, see \cite{gunzburger2014stochastic}.
Substituting Eq. \eqref{dataProjection} in the expression of $A_{j,j'}(\vec{y})$ yields, for all $j, j' = 1, \dots, J_h$,
\begin{align} \label{A_q}
A_{j,j'}(\vec{y}) &= \sum_{ \vec{q}\in \mathcal{J}(q)} \psi_{\vec{q}}(\vec{y})\int_D a_{\vec{q}}(\vec{x}) 
		      \nabla \phi_j(\vec{x}) \cdot \nabla \phi_{j'}(\vec{x}) d\vec{x} = \sum_{\vec{q}\in \mathcal{J}(q)} \psi_{\vec{q}}(\vec{y}) [A_{\vec{q}}]_{j,j'},
\end{align}
where $[A_{\vec{q}}]_{j,j'} = \int_D a_{\vec{q}}(\vec{x}) \nabla \phi_j(\vec{x}) \cdot \nabla \phi_{j'}(\vec{x}) d\vec{x}$
can be computed component-wise using a quadrature rule over $J_h$ elements. 
Replacing $A(\vec{y})$ in Eq. \eqref{linear_algebraic} with Eq. \eqref{A_q}, we obtain for all $\vec{p'} \in \mathcal{J}(p)$,  
\begin{equation} \label{linear_algebraic_2}
\sum_{\vec{p'}\in \mathcal{J}(p)} \sum_{\vec{q} \in \mathcal{J}(q)} 
\bigg[\int_{\Gamma} [A_{\vec{q}}] \psi_{\vec{q}}(\vec{y}) \psi_{\vec{p'}} (\vec{y}) 
\psi_{\vec{p}}(\vec{y}) \rho(\vec{y}) d\vec{y} \bigg] \vec{u}_{\vec{p}} = F_{\vec{p'}}.
\end{equation}
where $F_{\vec{p'}} = \int_{\Gamma} \vec{f} \psi_{\vec{p'}}(\vec{y}) \rho(\vec{y}) d \vec{y}.$
Let us define 
\begin{equation}
[G_{\vec{q}}]_{\vec{p'},\vec{p}} = \int_{\Gamma} \psi_{\vec{q}}(\vec{y}) \psi_{\vec{p'}} (\vec{y}) \psi_{\vec{p}}(\vec{y}) \rho(\vec{y}) d\vec{y}
\quad \mbox{and} \quad
K = \sum_{\vec{q} \in \mathcal{J}(q)} [G_{\vec{q}}] \bigotimes [A_{\vec{q}}],
\end{equation}
where $\bigotimes$ denotes the Kronecker product. Then Eq. \eqref{linear_algebraic_2} can be rewritten in matrix form as
\begin{equation} \label{stochastic_matrix}
 K \vec{u} = F.
\end{equation}
Note that the stochastic stiffness matrix $K$ in Eq. \eqref{stochastic_matrix} consists of ${(M_p)}^2$ block matrices, i.e.,
\begin{equation} \label{K_matrix}
K = 
\begin{bmatrix}
K_{1,1} & K_{1,2} & \cdots & K_{1, M_p} \\
K_{2,1} & K_{2,1} & \cdots & K_{2, M_p} \\
\vdots  & \vdots  & \ddots & \vdots  \\
K_{M_p,1} & K_{M_p,2} & \cdots & K_{M_p, M_p}
\end{bmatrix}
\end{equation}
where each block $K_{i,j}$ has the size of $[A_q]$, i.e. $J_h \times J_h$.

\section{Preconditioned GMRES solvers}\label{preconds}

\begin{figure}[t]
	\centering
	\includegraphics[width=1.0\textwidth]{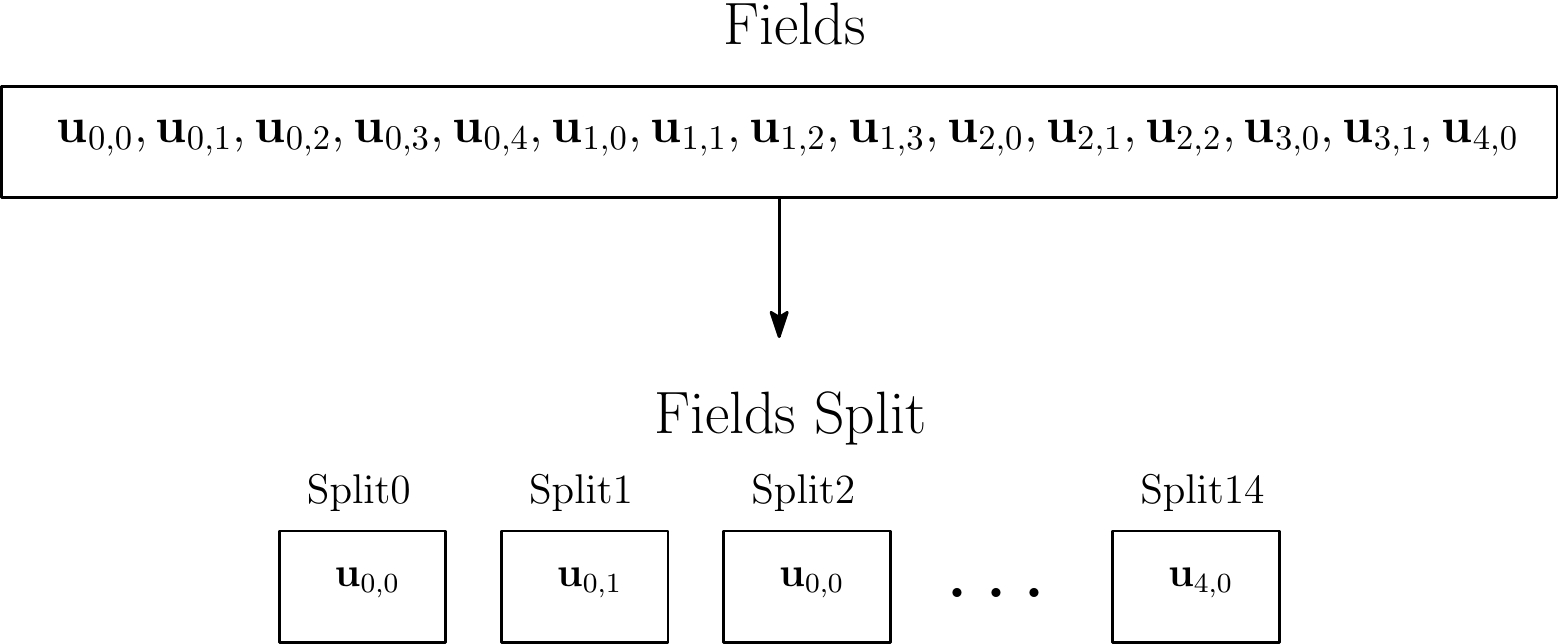}
	\caption{Schematics of the grouping procedure for the G-M-F-I preconditioned solver.}
	\label{fig:group_scheme_1}
\end{figure}

The system \eqref{stochastic_matrix} resulting from the discretization of the stochastic diffusion equation is solved with the GMRES algorithm
preconditioned with MG. On each multigrid level, the iterative solver below is employed as a smoother
\begin{align}\label{level_solver}
 \mathbf{u}^{(i+1)} = \mathbf{u}^{(i)} + P^{-1}(F-K\,\mathbf{u}^{(i)}),
\end{align}
where $i=1,\ldots,i_{\max}$, $i_{\max}=1$ unless otherwise stated, and $P^{-1}$ represents the action of an appropriate preconditioner.
At the coarsest level, a direct solver is used in place of Eq. \eqref{level_solver}. 
Two different strategies are used to select $P$. 
With the first, $P$ is the incomplete LU (ILU) factorization of the whole matrix $K$, whereas with the second $P$ is a matrix obtained according to an FS strategy, as explained next. With this second approach, it is required to invert diagonal blocks of a matrix, and these inverse matrices are approximated using an ILU factorization.
The solver obtained using the ILU factorization directly on $K$ is labeled G-M-I (because GMRES is preconditioned by MG with level solver in Eq. \eqref{level_solver} preconditioned by ILU), whereas the one obtained using FS is denoted by G-M-F-I (because GMRES is preconditioned by MG with level solver in Eq. \eqref{level_solver} preconditioned by FS and then ILU is used to approximate the inverses of the diagonal blocks). 
Considering that the G-M-I approach is fairly standard and well established, we focus on the field-split type preconditioners, which are the main novelty of this work. The G-M-F-I is obtained
with the most straightforward splitting strategy: to each stochastic mode (also referred to as {\textit{field}}) we associate one split.
For instance, if we have 15 stochastic fields, i.e. $M_p = 15$,  
then 15 splits are generated, one for each mode.
Hence, referring to the stochastic stiffness matrix $K$ in \eqref{K_matrix}, with G-M-F-I the iterative solver in \eqref{level_solver} is preconditioned with the matrix $P$ below
 \begin{equation}  \label{P_preconditioner}
P= 
\begin{bmatrix}
K_{1,1} & 0 & \cdots & 0 \\
K_{2,1} & K_{2,1} &\cdots & 0 \\
\vdots  & \vdots  & \ddots & \vdots  \\
K_{M_p,1} & K_{M_p,2} & \cdots & K_{M_p, M_p}
\end{bmatrix}.
\end{equation}
This choice of $P$ is usually refereed to as {\it{left multiplicative field-split preconditioner}}. 
Because $P$ is a lower triangular matrix, computing its inverse
requires to invert only the diagonal blocks $K_{i,i}$ for  $i =0, \dots,M_p$. 
As mentioned above, each inverse is approximated using an ILU factorization.
Schematics of the field-splits that make up the G-M-F-I are shown in Figure \ref{fig:group_scheme_1}.
A preliminary numerical investigation showed that the fields obtained from the splitting just described can be grouped back together to further enhance the performances of the preconditioner.
Namely, a robust and efficient preconditioner can be obtained grouping the fields according to the leading digits of the $N$-tuples of $\mathcal{J}(p)$. The preconditioned solver that takes advantage of this approach is labeled G-M-Fg-I.
With this new strategy, we aim to reduce the number of splits for the matrix $K$, which otherwise grows fast with $N$ and $p$, i.e. $ M_p = (N+p)!/(N! p!)$, so that the resulting G-M-Fg-I algorithm is faster and more robust than G-M-F-I.

\begin{figure}[t]
	\centering
	\includegraphics[width=1.0\textwidth]{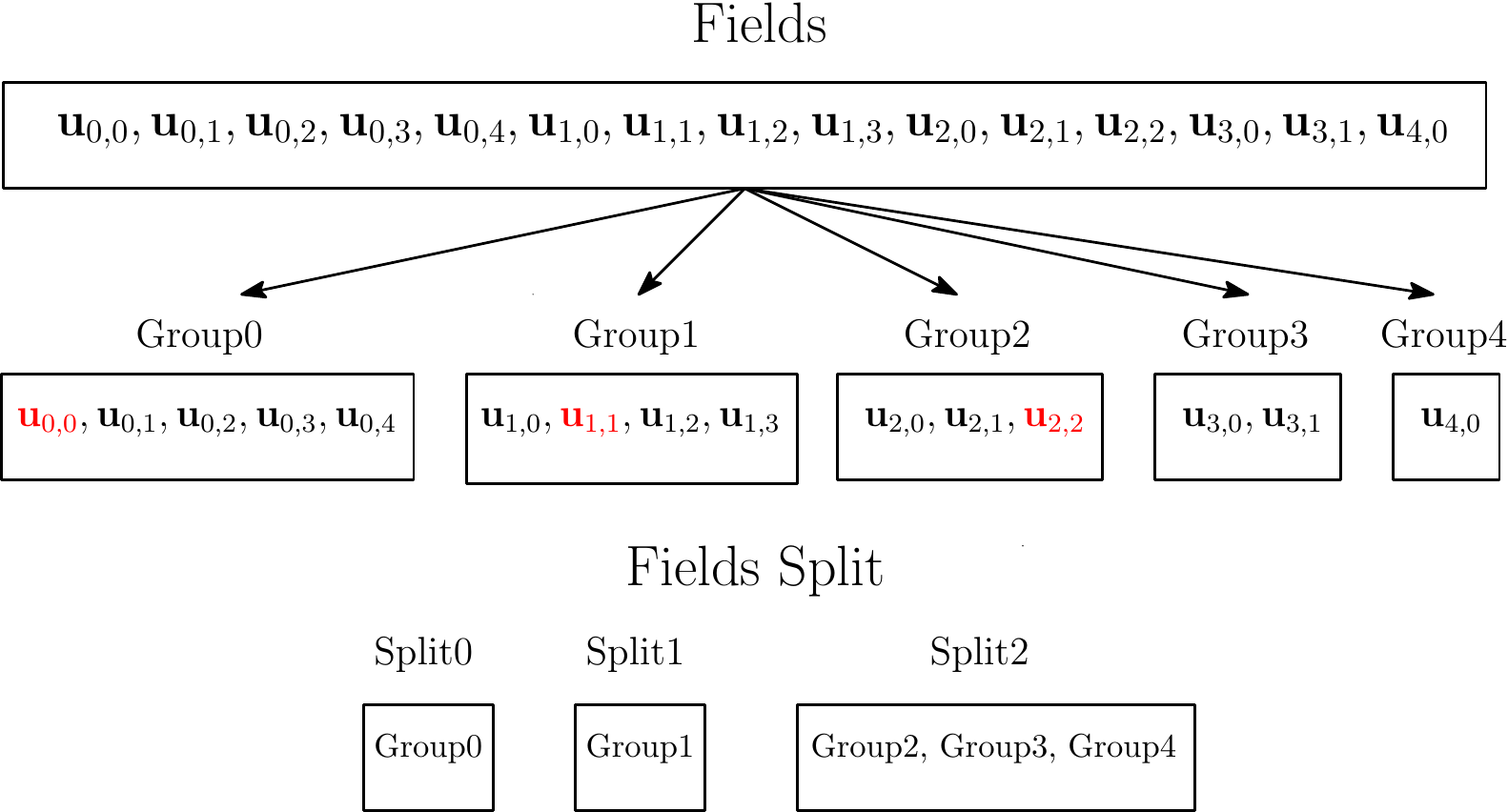}
	\caption{Schematics of the grouping procedure for the G-M-Fg-I preconditioned solver.}
	\label{fig:group_scheme}
\end{figure}

To illustrate this alternative grouping procedure, we consider the case of $N=2$ and $p = 4$, but the technique can be easily generalized to other values of $N$ and $p$. 
For the values of $N$ and $p$ chosen, we have $M_p = 15$, hence a total of 15 fields. 
Next we group the fields based on the same first digit in the multi-index and generate 5 groups: 
Group0, Group1, Group2, Group3, and Group4, as showed in Figure \ref{fig:group_scheme}.

For any group, if there exists one field whose multi-index has all identical digits, i.e. the red fields in Figure \ref{fig:group_scheme},  
this group will form an new split, otherwise it will be included in the preceding split.
With this principle, three splits are created:
Split1 consisting of Group0, Split2 consisting of Group1, and 
Split3 consisting of Group2, Group3, and Group4.
With this grouping method, the preconditioner $P$ used in Eq. \eqref{level_solver} for G-M-Fg-I can be written as
\begin{equation} \label{P_group}
P= 
\begin{bmatrix}
K_{g0,g0} & 0 & 0 \\
K_{g1,g0} & K_{g1,g1} & 0  \\
K_{g2,g0} & K_{g2,g1} & K_{g2,g2} \\
\end{bmatrix}.
\end{equation}
where the inverse of the diagonal sub-blocks $K_{gi,gi}$, $i=1,\dots,3$, are again approximated by the inverse of their ILU decomposition.

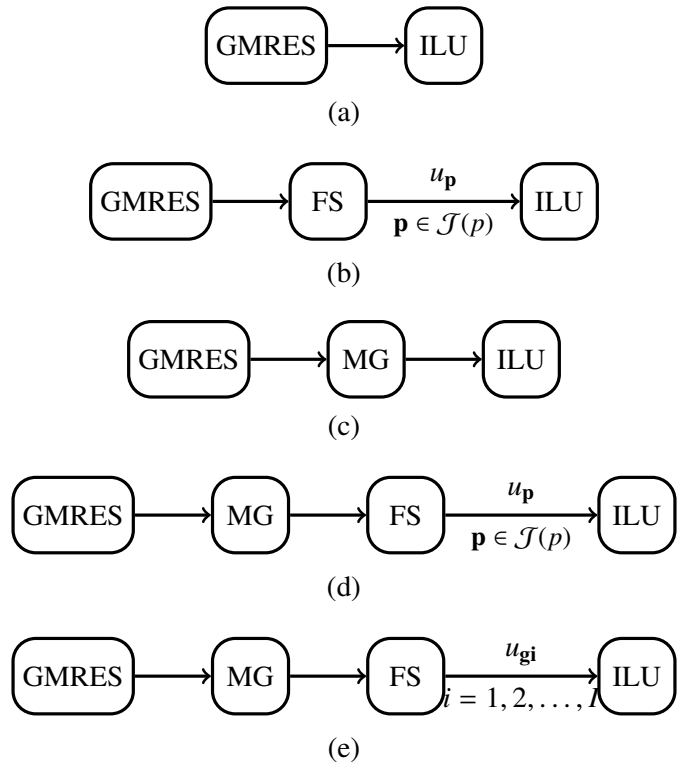
\begin{figure}[!htbp]
\begin{minipage}{1.0\textwidth}
\begin{center}
 
\begin{minipage}{1.0\textwidth}
\begin{center}
\begin{tikzpicture}[
node distance=3cm,
     terminal/.style={
      rectangle,rounded corners=3mm,minimum size=1cm,
       very thick,draw=black!100!blue,
       top color=white,
       }
]
{\node (gmres) [terminal] {GMRES}; }
{\node (ilu)   [terminal,position=0:{1.cm} from gmres] { ILU }; }
{\path (gmres) edge[->,very thick] (ilu); }
\end{tikzpicture}\\
(a)
\end{center}
\end{minipage}

\vspace{0.4cm}

\begin{minipage}{1.0\textwidth}
\begin{center}
\begin{tikzpicture}[
     node distance = 3cm,
     terminal/.style={
      rectangle,rounded corners=3mm,minimum size=1cm,
       very thick,draw=black!100!blue,
       top color=white,
       }
]
{\node (gmres) [terminal] {GMRES}; }
{\node (fs) [terminal,position=0:{1.cm} from gmres] { FS }; }
{\node (ilu1)  [terminal,position= 0:{\nodeDist} from fs] { ILU }; }

{\path (gmres) edge[->,very thick] (fs); }
{\path (fs) edge[->,very thick] node [midway, above, sloped] (TextNode2) {$u_{\mathbf{p}}$} node [midway, below, sloped] (TextNode2) {\small $\mathbf{p} \in \mathcal{J}(p)$} (ilu1); }
\end{tikzpicture} \\
(b)
\end{center}
\end{minipage}

\vspace{0.4cm}

\begin{minipage}{1.0\textwidth}
\begin{center}
\begin{tikzpicture}[
node distance=3cm,
     terminal/.style={
      rectangle,rounded corners=3mm,minimum size=1cm,
       very thick,draw=black!100!blue,
       top color=white,
       }
]
{\node (gmres) [terminal] {GMRES}; }
{\node (gmg) [terminal,position=0:{1.cm} from gmres] { MG }; }
{\node (ilu)  [terminal,position=0:{1.cm} from gmg] { ILU }; }
{\path (gmres) edge[->,very thick] (gmg); }
{\path (gmg) edge[->,very thick] (ilu); }
\end{tikzpicture}\\
(c)
\end{center}
\end{minipage}

\vspace{0.4cm}

\begin{minipage}{1.0\textwidth}
\begin{center}
\begin{tikzpicture}[
     node distance = 3cm,
     terminal/.style={
      rectangle,rounded corners=3mm,minimum size=1cm,
       very thick,draw=black!100!blue,
       top color=white,
       }
]
{\node (gmres) [terminal] {GMRES}; }
{\node (gmg) [terminal,position=0:{1.cm} from gmres] { MG}; }
{\node (fs) [terminal,position=0:{1.cm} from gmg] { FS }; }
{\node (ilu1)  [terminal,position= 0:{\nodeDist} from fs] { ILU }; }

{\path (gmres) edge[->,very thick] (gmg); }
{\path (gmg) edge[->,very thick] (fs); }
{\path (fs) edge[->,very thick] node [midway, above, sloped] (TextNode2) {$u_{\mathbf{p}}$} node [midway, below, sloped] (TextNode2) {\small $\mathbf{p} \in \mathcal{J}(p)$} (ilu1); }
\end{tikzpicture} \\
(d)
\end{center}
\end{minipage}

\vspace{0.4cm}

\begin{minipage}{1.0\textwidth}
\begin{center}
\begin{tikzpicture}[
     node distance = 3cm,
     terminal/.style={
      rectangle,rounded corners=3mm,minimum size=1cm,
       very thick,draw=black!100!blue,
       top color=white,
       }
]
{\node (gmres) [terminal] {GMRES}; }
{\node (gmg) [terminal,position=0:{1.cm} from gmres] {  MG}; }
{\node (fs) [terminal,position=0:{1.cm} from gmg] { FS }; }
{\node (ilu1)  [terminal,position= 0:{\nodeDist} from fs] { ILU }; }

{\path (gmres) edge[->,very thick] (gmg); }
{\path (gmg) edge[->,very thick] (fs); }
{\path (fs) edge[->,very thick] node [midway, above, sloped] (TextNode2) {$u_{\mathbf{gi}}$} node [midway, below, sloped] (TextNode2) {$i = 1, 2,\dots, I$} (ilu1); }
\end{tikzpicture} \\
(e)
\end{center}
\end{minipage}
\end{center}
\end{minipage}
\caption{Building blocks that make up each preconditioned solver considered in the numerical investigation: 
(a) G-I 
(b) G-F-I 
(c) G-M-I 
(d) G-M-F-I 
(e) G-M-Fg-I. 
 } 
\label{fig:preconditioners_chart}
\end{figure}

To have a well assorted set of preconditioners to compare, 
we also apply ILU  and FS+ILU directly to the GMRES solver, producing the methods here denoted as G-I and G-F-I, respectively.
These preconditioners can be immediately derived from the ones already described considering only one level for MG without any coarse grid correction. 
In order to boost the performance of the G-F-I preconditioned solver, and only for this case, we consider $i_{\max}=50$ in Eq. \eqref{level_solver}
to obtain a more accurate action of the inverse of each diagonal block $K_{i,i}$. 
We did not adopt this choice for G-M-F-I and G-M-Fg-I, because none of them benefited from it in terms of computational time.
The alternative G-Fg-I grouping is not considered because in a preliminary numerical investigation it exhibited inferior convergence properties than the G-F-I. A summary of the methods discussed is shown in Figure \ref{fig:preconditioners_chart}.

\section{Computational Investigation}\label{numres}
In this section, we apply the five preconditioned solvers discussed in the previous section to solve the two dimensional stochastic diffusion problem in \eqref{Strong_PS} with $f\equiv1$.
The physical domain is a unit square with Dirichlet zero boundary conditions.  
For the approximation of the stochastic lognormal coefficient $a(\vec{x},\vec{y})$ with the truncated KL expansion $a_{KL}(\vec{x},\vec{y})$ in Eq. \eqref{KL_gamma} we choose $\mu_{\gamma} = 0$, $a_{min} = 0.01$, and $L_c=0.1$. 
All the numerical simulations have been carried out on a Dell OptiPlex 760 with SFF/Core 2 Duo E8400 @ 3.00 GHz 
using FEMuS, a in-house open-source finite element C++ library built on top of PETSc \cite{balay2012petsc}. 
In the GMRES solver, the absolute error and relative error for the scaled preconditioned residual of each variable 
are set to $1.0 \times 10 ^{-10}$ and $1.0 \times 10 ^{-10}$, respectively. 
A V-cycle is used for the multigrid, with one pre and post-smoothing iteration, i.e. $i_{\max}=1$ in Eq. \eqref{level_solver}, except for the G-F-I for which $i_{\max}=50$.
Note that our algorithm is not limited to the V-cycle multigrid but also applicable to the W-cycle or the F-cycle. 

We first investigate the effect of mesh refinement on the performances of the preconditioners.
For this case, we fix the standard deviation $\sigma_{\gamma} = 0.08 $ in \eqref{standard_deviation},  
$p=4$ in \eqref{multivariate_polynomial}, $q=5$ in \eqref{dataProjection}, and the stochastic dimension $ N = 2$ in \eqref{aKL}. 
Starting with a coarse mesh with $2 \times 2$ quadrilateral elements indicated as level $L= 1$, 
we perform a midpoint refinement procedure to obtain the subsequent levels. For instance, $L=3$ refers to a mesh with $8 \times 8$ cells obtained refining $L-1$ times the original coarse mesh.
The numerical results for the computational time and the number of GMRES iterations are reported 
in Table \ref{mesh_variation}.
The number of GMRES iterations increases dramatically for G-I and noticeably for G-F-I, showing 
a rate of convergence that depends on $L$.
On the other hand, the number of GMRES iterations remains stable for the MG-type preconditioned solvers, namely
G-M-I, G-M-F-I, and G-M-Fg-I, displaying a rate of convergence that is independent of the mesh size. 
G-M-F-I and G-M-Fg-I are the ones that present the best results in terms of computational time among the methodologies tested,
being at least twice faster than the other methods. 

\begin {table}[!t]
\begin{center}
	\begin{tabular}{|c|c|c|c|c|c|c|c|} \hline	
	Precond. & {G-I} & {G-F-I} & {G-M-I} & {G-M-F-I} & {G-M-Fg-I} \\ \hline 
	     L  & \shortstack{\\CPU Time [s] \\ (GMRES iter.)}  & \shortstack{\\CPU Time [s] \\ (GMRES iter.)}  & \shortstack{\\CPU Time [s] \\ (GMRES iter.)}  & \shortstack{CPU Time [s] \\ (GMRES iter.)} & \shortstack{CPU Time [s] \\ (GMRES iter.)} \\ \hline   
	     3 & \shortstack{\\0.27   \\(31)}   & \shortstack{\\0.09 \\ (3)}  & \shortstack{\\0.27  \\ (7)}  & \shortstack{\\0.13   \\ (7)} & \shortstack{\\0.13   \\ (7)} \\ \hline 
	     4 & \shortstack{\\1.66   \\(60)}   & \shortstack{\\0.61 \\ (7)}  & \shortstack{\\1.28  \\ (7)}  & \shortstack{\\0.48   \\ (7)} & \shortstack{\\0.48   \\ (7)} \\ \hline 
	     5 & \shortstack{\\11.60  \\(129)}  & \shortstack{\\4.74 \\ (16)} & \shortstack{\\5.65  \\ (7)}  & \shortstack{\\1.84   \\ (6)} & \shortstack{\\1.83   \\ (6)} \\ \hline 
	     6 & \shortstack{\\142.97 \\(363)}  & \shortstack{\\59.84\\ (40)} & \shortstack{\\26.81 \\ (7)}  & \shortstack{\\10.20  \\ (6)} & \shortstack{\\10.08   \\ (6)} \\ \hline 
	\end{tabular}
\end{center}
\caption{CPU times and GMRES iterations as the size of the mesh is decreased, considering  $\sigma_{\gamma} = 0.08$, $p=4$, $q=5$, and $N = 2$.}
\label{mesh_variation}
\end{table}

\begin {table}[!htpb] 
\begin{center}
	\begin{tabular}{|c|c|c|c|c|c|c|c|} \hline	
	 Precond. & {G-I} & {G-F-I} & {G-M-I} & {G-M-F-I} & {G-M-Fg-I} \\ \hline 
	     ~$\sigma_{\gamma}$~ & \shortstack{\\CPU Time [s] \\ (GMRES iter.)} & \shortstack{\\CPU Time [s] \\ (GMRES iter.)}  & \shortstack{\\CPU Time [s] \\ (GMRES iter.)}  & \shortstack{CPU Time [s] \\ (GMRES iter.)} & \shortstack{CPU Time [s] \\ (GMRES iter.)} \\ \hline  
	     0.08 & \shortstack{\\25.92 \\ (129)} & \shortstack{\\7.18 \\(16)}  & \shortstack{13.85 \\(7)} & \shortstack{\\3.81 \\ (6)}  & \shortstack{\\3.81 \\ (6)} \\ \hline
	     0.1  & \shortstack{\\26.92 \\ (135)} & \shortstack{\\7.70 \\(17)}  & \shortstack{13.86 \\(7)} & \shortstack{\\3.82 \\ (6)}  & \shortstack{\\3.83 \\ (6)} \\ \hline 
	     0.2  & \shortstack{\\28.25 \\ (145)} & \shortstack{\\8.82 \\(19)}  & \shortstack{13.86 \\(7)} & \shortstack{\\3.82 \\ (6)}  & \shortstack{\\3.83 \\ (6)} \\ \hline 
	     0.4  & \shortstack{\\30.86 \\ (163)} & \shortstack{\\9.48 \\(20)}  & \shortstack{13.87 \\(7)} & \shortstack{\\4.24 \\ (7)}  & \shortstack{\\4.29 \\ (7)} \\ \hline 
	     0.8  & \shortstack{\\35.02 \\ (195)} & \shortstack{\\11.59 \\(25)} & \shortstack{13.88 \\(7)} & \shortstack{\\4.66 \\ (8)}  & \shortstack{\\4.68 \\ (8)} \\ \hline 
	     1.0  & \shortstack{\\36.97 \\ (207)} & \shortstack{\\12.82 \\(28)} & \shortstack{13.88 \\(7)} & \shortstack{\\5.53 \\ (10)} & \shortstack{\\4.72 \\ (8)} \\ \hline 
	     1.4  & \shortstack{\\40.07 \\ (232)} & \shortstack{\\14.22 \\(31)} & \shortstack{13.89 \\(7)} & \shortstack{\\8.42 \\ (17)} & \shortstack{\\5.11 \\ (9)} \\ \hline
	     1.6  & \shortstack{\\42.48 \\ (246)} & \shortstack{\\16.71 \\(37)} & \shortstack{14.30 \\(8)} & \shortstack{\\11.10\\ (23)} & \shortstack{\\5.14 \\ (9)} \\ \hline
	\end{tabular}
\end{center}
\caption{CPU times and GMRES iterations as the input standard deviation $\sigma_{\gamma}$ is increased, considering $L=5$, $p=q=5$, and $N = 2$}
\label{sigma_variation}
\end{table}
\begin{figure}[!htpb] 
\begin{minipage}{1.0\textwidth}
\begin{center}
\begin{minipage}{0.49\textwidth}
\begin{center}
\includegraphics[width=1.0\linewidth]{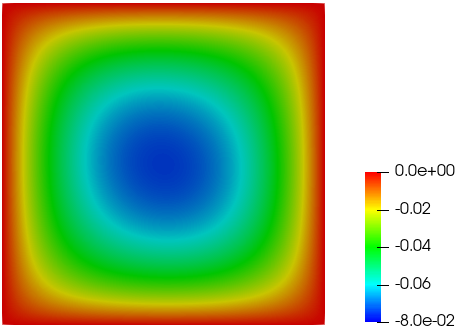}\\
$u_{0,0}$
\end{center}
\end{minipage}
\begin{minipage}{0.49\textwidth}
\begin{center}
\includegraphics[width=1.0\linewidth]{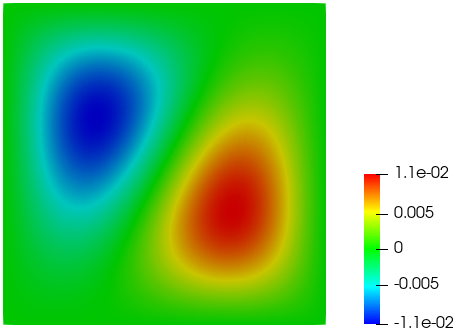} \\
$u_{0,1}$
\end{center}
\end{minipage}
\end{center}
\end{minipage}
\caption{Graphs of $u_{0,0}$ and $u_{0,1}$ for the case $\sigma_{\gamma} = 1.6$, $L=5$, $p=q=5$, and $N = 2$. } 
\label{usg_solution}
\end{figure}
Next, the effect of the input standard deviation $\sigma_{\gamma}$ on the robustness and efficiency of the 
preconditioned solvers is studied.
The input standard deviation is increased from $0.08$ to $1.6$, 
while fixing $L=5$, $p=q=5$, and the stochastic dimension $N=2$.  
The finite element solutions associated with the first and the second stochastic mode, i.e. ${u}_{0,0}$ and ${u}_{0,1}$, obtained with G-M-F-I are shown in Figure \ref{usg_solution} for the case of $\sigma_{\gamma} = 1.6$. Recall that $u_{\vec{p}}$ is defined in Eq. \eqref{u_p}.
Results are reported in Table \ref{sigma_variation}: concerning the number of GMRES iterations, we see that for G-I and G-F-I
there is a monotone increase with $\sigma_{\gamma}$, 
suggesting a rate of convergence that depends on this parameter.
\begin {table}[!t] 
\begin{center}
	\begin{tabular}{|c|c|c|c|c|c|c|c|} \hline	
	 Precond. & {G-I} & {G-F-I} & {G-M-I} & {G-M-F-I} & {G-M-Fg-I} \\ \hline 
	     ~$p$~ & \shortstack{\\CPU Time [s] \\ (GMRES iter.)}  & \shortstack{\\CPU Time [s] \\ (GMRES iter.)}  & \shortstack{\\CPU Time [s] \\ (GMRES iter.)}  & \shortstack{CPU Time [s] \\ (GMRES iter.)} & \shortstack{CPU Time [s] \\ (GMRES iter.)} \\ \hline 	     
	     4 & \shortstack{\\11.60 \\ (129)} & \shortstack{\\4.74  \\(16)} & \shortstack{\\ 5.65 \\ (7)} & \shortstack{\\1.84  \\(6)} & \shortstack{\\1.83  \\(6)} \\ \hline 
	     5 & \shortstack{\\25.92 \\ (129)} & \shortstack{\\7.18  \\(16)} & \shortstack{\\13.85 \\ (7)} & \shortstack{\\3.81  \\(6)} & \shortstack{\\3.81  \\(6)} \\ \hline 
	     6 & \shortstack{\\52.50 \\ (129)} & \shortstack{\\10.33 \\(16)} & \shortstack{\\30.05 \\ (7)} & \shortstack{\\7.01  \\(6)} & \shortstack{\\6.89  \\(6)} \\ \hline 
	     7 & \shortstack{\\97.90 \\ (129)} & \shortstack{\\14.47 \\(16)} & \shortstack{\\61.51 \\ (7)} & \shortstack{\\13.37 \\(6)} & \shortstack{\\14.70 \\(6)} \\ \hline 
	\end{tabular}
\end{center}
\caption{CPU times and GMRES iterations with $\sigma_{\gamma} = 0.08 $ as the integer $p$ in Eq. \eqref{multivariate_polynomial} is increased, considering level = 5,  $q = 5$, and $N = 2$.}
\label{p_var_low}
\end{table}

\begin {table}[!t] 
\begin{center}
	\begin{tabular}{|c|c|c|c|c|c|c|c|} \hline	
	 Precond. & {G-I} & {G-F-I} & {G-M-I} & {G-M-F-I} & {G-M-Fg-I} \\ \hline
	     ~$p$~ & \shortstack{\\CPU Time [s] \\ (GMRES iter.)}  & \shortstack{\\CPU Time [s] \\ (GMRES iter.)}  & \shortstack{\\CPU Time [s] \\ (GMRES iter.)}  & \shortstack{CPU Time [s] \\ (GMRES iter.)} & \shortstack{CPU Time [s] \\ (GMRES iter.)} \\ \hline 	     
	     4 & \shortstack{\\15.95  \\ (193)} & \shortstack{\\7.00  \\(24)} & \shortstack{\\ 5.66 \\ (7)} & \shortstack{\\2.24  \\(8)}  & \shortstack{\\2.20  \\(8)} \\ \hline 
	     5 & \shortstack{\\35.02  \\ (195)} & \shortstack{\\11.59 \\(25)} & \shortstack{\\13.88 \\ (7)} & \shortstack{\\4.66  \\(8)}  & \shortstack{\\4.68  \\(8)} \\ \hline 
	     6 & \shortstack{\\69.95  \\ (195)} & \shortstack{\\18.10 \\(26)} & \shortstack{\\30.05 \\ (7)} & \shortstack{\\8.62  \\(8)}  & \shortstack{\\8.44  \\(8)} \\ \hline 
	     7 & \shortstack{\\128.19 \\ (195)} & \shortstack{\\27.38 \\(27)} & \shortstack{\\61.52 \\ (7)} & \shortstack{\\17.34 \\(8)}  & \shortstack{\\17.26 \\(8)} \\ \hline 
	\end{tabular}
\end{center}
\caption{CPU times and GMRES iterations with $\sigma_{\gamma} = 0.8 $ as the integer $p$ in Eq. \eqref{multivariate_polynomial} is increased, considering level = 5,  $q = 5$, and $N = 2$.}
\label{p_var_middle}
\end{table}

\begin {table}[!t] 
\begin{center}
	\begin{tabular}{|c|c|c|c|c|c|c|c|} \hline	
	 Precond. & {G-I} & {G-F-I} & {G-M-I} & {G-M-F-I} & {G-M-Fg-I} \\ \hline
	     ~$p $~ & \shortstack{\\CPU Time [s] \\ (GMRES iter.)}  & \shortstack{\\CPU Time [s] \\ (GMRES iter.)}  & \shortstack{\\CPU Time [s] \\ (GMRES iter.)}  & \shortstack{CPU Time [s] \\ (GMRES iter.)} & \shortstack{CPU Time [s] \\ (GMRES iter.)}\\ \hline 	     
	     4 & \shortstack{\\18.32  \\ (229)} & \shortstack{\\8.35  \\(29)} & \shortstack{\\5.66  \\ (7)} & \shortstack{\\3.46  \\(14)} & \shortstack{\\2.41  \\(9)} \\ \hline 
	     5 & \shortstack{\\40.07  \\ (232)} & \shortstack{\\14.22 \\(31)} & \shortstack{\\13.89 \\ (7)} & \shortstack{\\8.42  \\(17)} & \shortstack{\\5.11  \\(9)} \\ \hline 
	     6 & \shortstack{\\80.18  \\ (234)} & \shortstack{\\24.04 \\(35)} & \shortstack{\\30.05 \\ (7)} & \shortstack{\\18.90 \\(21)} & \shortstack{\\9.99  \\(10)}\\ \hline 
	     7 & \shortstack{\\145.64 \\ (235)} & \shortstack{\\46.37 \\(47)} & \shortstack{\\61.52 \\ (7)} & \shortstack{\\59.41 \\(33)} & \shortstack{\\22.31 \\(10)}\\ \hline 
	\end{tabular}
\end{center}
\caption{CPU times and GMRES iterations with $\sigma_{\gamma} = 1.4 $ as the integer $p$ in Eq. \eqref{multivariate_polynomial} is increased, considering level = 5,  $q = 5$, and $N = 2$.}
\label{p_var_high}
\end{table}
More interestingly, for G-M-F-I 
the number of GMRES iteration remains stable for $\sigma_{\gamma} \leq  0.8$, but it grows consistently for $\sigma_{\gamma} > 0.8$, going from 6 iterations with $\sigma_{\gamma} =  0.08$ to 23 iterations with $\sigma_{\gamma} =  1.6$.
A similar behavior for a preconditioned GMRES solver was observed in \cite{ullmann2012efficient}.
On the other hand, the G-M-Fg-I is practically insensitive to variations of  $\sigma_{\gamma}$ going from 6 to 9 iterations for the values of $\sigma_{\gamma}$ considered. 
This behavior shows that the grouping strategy adopted for G-M-Fg-I has significantly improved the robustness of the solver over the straightforward grouping used in G-M-F-I.
The G-M-I is also insensitive to variations of $\sigma_{\gamma}$.
Concerning the CPU times, 
the G-M-F-I and G-M-Fg-I methods are generally faster than the others. 
For $\sigma_{\gamma} \leq 0.8 $, the performance of G-M-F-I and G-M-Fg-I are comparable, 
whereas for $\sigma_{\gamma} > 0.8 $, G-M-Fg-I outperforms G-M-F-I.
Such a behavior can be attributed to the difference in number of GMRES iterations between the two methods.
For instance, for $\sigma_{\gamma} = 1.6$, 
the G-M-Fg-I is about twice as fast as G-M-F-I,
because G-M-Fg-I has half the iterations that G-M-F-I has.
The G-M-Fg-I is always the fastest method among those considered.
The CPU time of G-M-I grows slowly, although it is three times as large as G-M-Fg-I.
The G-F-I is faster than the G-M-I for small $\sigma_{\gamma}$, although it becomes slower large $\sigma_{\gamma}$.
The CPU time of G-I is always the worse among the methods considered.
To summarize, regarding the dependence on $\sigma_{\gamma}$, G-M-Fg-I is the fastest method and is optimal, G-M-F-I is the second fastest but is not optimal, G-M-I is optimal but slower than the previous two. G-I and G-F-I deteriorate both in terms of iterations and CPU time.

We continue our analysis investigating the effect of the dimension of the stochastic finite dimensional space varying the value of $p$.
We vary $p$  from $4$ to $7$ and set $L=5$,  $q=5$, and the stochastic input dimension $N = 2$. 
Three specific values of $\sigma_{\gamma}$ are selected, namely $0.08$, $0.8$, and $1.4$, 
that represent a relatively small, medium, and large standard deviation. 
Results are reported in Tables \ref{p_var_low}, 
\ref{p_var_middle}, and \ref{p_var_high}, respectively.  
For the small and medium $\sigma_{\gamma}$ in Tables \ref{p_var_low} and \ref{p_var_middle} respectively, 
we observe that the number of GMRES iterations remains stable for all methods as $p$ increases, 
meaning that the rate of convergence is independent of $p$. 
The MG-type methods always perform better than G-I and G-F-I 
in terms of number of iterations. 
Concerning the CPU time, G-M-F-I and G-M-Fg-I are again the most efficient, moreover G-F-I outperforms G-M-I.
A different situation happens for the large value of $\sigma_{\gamma}$ in Table \ref{p_var_high}, 
where the number of GMRES iterations increases significantly for G-F-I and G-M-F-I  
whereas for G-I, G-M-I and G-M-Fg-I it remains practically stable. 
It is also found that G-M-Fg-I outperforms all other methods in terms of computational time.
For instance if $p = 7$, G-M-Fg-I is at least 2 times faster than G-F-I, G-M-I, and G-M-F-I and 6.5 times faster than G-I.
\begin {table}[!ht] 
\begin{center}
	\begin{tabular}{|c|c|c|c|c|c|c|c|} \hline	
	 Precond. & {G-I} & {G-F-I} & {G-M-I} & {G-M-F-I} & {G-M-Fg-I}\\ \hline
	     ~$q$~ & \shortstack{\\CPU Time [s] \\ (GMRES iter.)}  & \shortstack{\\CPU Time [s] \\ (GMRES iter.)}  & \shortstack{\\CPU Time [s] \\ (GMRES iter.)}  & \shortstack{CPU Time [s] \\ (GMRES iter.)} & \shortstack{CPU Time [s] \\ (GMRES iter.)}  \\ \hline 	     
	     4 & \shortstack{\\25.39  \\ (129)} & \shortstack{\\7.11 \\(16)} & \shortstack{\\13.18 \\ (7)} & \shortstack{\\3.78 \\(6)}  & \shortstack{\\3.81 \\(6)} \\ \hline 
	     5 & \shortstack{\\25.92  \\ (129)} & \shortstack{\\7.18 \\(16)} & \shortstack{\\13.85 \\ (7)} & \shortstack{\\3.81 \\(6)}  & \shortstack{\\3.81 \\(6)} \\ \hline 
	     6 & \shortstack{\\25.92  \\ (129)} & \shortstack{\\7.18 \\(16)} & \shortstack{\\13.85 \\ (7)} & \shortstack{\\3.81 \\(6)}  & \shortstack{\\3.82 \\(6)} \\ \hline 
	     7 & \shortstack{\\25.93  \\ (129)} & \shortstack{\\7.18 \\(16)} & \shortstack{\\13.91 \\ (7)} & \shortstack{\\3.81 \\(6)}  & \shortstack{\\3.82 \\(6)} \\ \hline 
	\end{tabular}
\end{center}
\caption{CPU times and GMRES iterations with $\sigma_{\gamma} = 0.08 $ as the integer $q$ in Eq. \eqref{dataProjection} is increased, with $L = 5$, $p = 5$, and $N = 2$.}
\label{q_var_low}
\end{table}

\begin {table}[!ht] 
\begin{center}
	\begin{tabular}{|c|c|c|c|c|c|c|c|} \hline	
	 Precond. & {G-I} & {G-F-I} & {G-M-I} & {G-M-F-I} & {G-M-Fg-I}\\ \hline
	     ~$q$~ & \shortstack{\\CPU Time [s] \\ (GMRES iter.)}  & \shortstack{\\CPU Time [s] \\ (GMRES iter.)}  & \shortstack{\\CPU Time [s] \\ (GMRES iter.)}  & \shortstack{CPU Time [s] \\ (GMRES iter.)} & \shortstack{CPU Time [s] \\ (GMRES iter.)}\\ \hline 	     
	     4 & \shortstack{\\34.53  \\ (195)} & \shortstack{\\11.49 \\(25)} & \shortstack{\\13.22 \\ (7)} & \shortstack{\\4.65 \\(8)} & \shortstack{\\4.68 \\(8)}\\ \hline 
	     5 & \shortstack{\\35.02  \\ (195)} & \shortstack{\\11.59 \\(25)} & \shortstack{\\13.88 \\ (7)} & \shortstack{\\4.66 \\(8)} & \shortstack{\\4.68 \\(8)}\\ \hline 
	     6 & \shortstack{\\35.02  \\ (195)} & \shortstack{\\11.60 \\(25)} & \shortstack{\\13.90 \\ (7)} & \shortstack{\\4.66 \\(8)} & \shortstack{\\4.70 \\(8)}\\ \hline 
	     7 & \shortstack{\\35.02  \\ (195)} & \shortstack{\\11.60 \\(25)} & \shortstack{\\13.90 \\ (7)} & \shortstack{\\4.66 \\(8)} & \shortstack{\\4.70 \\(8)}\\ \hline 
	\end{tabular}
\end{center}
\caption{CPU times and GMRES iterations with $\sigma_{\gamma} = 0.8 $ as the integer $q$ in Eq. \eqref{dataProjection} is increased, $L = 5$, $p = 5$, and $N = 2$.}
\label{q_var_middle}
\end{table}

\begin {table}[!ht] 
\begin{center}
	\begin{tabular}{|c|c|c|c|c|c|c|c|} \hline	
	 Precond. & {G-I} & {G-F-I} & {G-M-I} & {G-M-F-I} & {G-M-Fg-I}\\ \hline
	     ~$q$~ & \shortstack{\\CPU Time [s] \\ (GMRES iter.)}  & \shortstack{\\CPU Time [s] \\ (GMRES iter.)}  & \shortstack{\\CPU Time [s] \\ (GMRES iter.)}  & \shortstack{CPU Time [s] \\ (GMRES iter.)} & \shortstack{CPU Time [s] \\ (GMRES iter.)}\\ \hline 	     
	     4 & \shortstack{\\39.65  \\ (232)} & \shortstack{\\14.09 \\(31)} & \shortstack{\\13.23 \\ (7)} & \shortstack{\\7.96 \\(16)}  & \shortstack{\\5.10 \\(9)} \\ \hline 
	     5 & \shortstack{\\40.07  \\ (232)} & \shortstack{\\14.22 \\(31)} & \shortstack{\\13.89 \\ (7)} & \shortstack{\\8.42 \\(17)}  & \shortstack{\\5.11  \\(9)} \\ \hline 
	     6 & \shortstack{\\40.07  \\ (232)} & \shortstack{\\14.22 \\(31)} & \shortstack{\\13.90 \\ (7)} & \shortstack{\\8.42 \\(17)}  & \shortstack{\\5.11  \\(9)} \\ \hline 
	     7 & \shortstack{\\40.15  \\ (232)} & \shortstack{\\14.23 \\(31)} & \shortstack{\\13.90 \\ (7)} & \shortstack{\\8.42 \\(17)}  & \shortstack{\\5.11  \\(9)} \\ \hline 
	\end{tabular}
	\end{center}
\caption{CPU times and GMRES iterations with $\sigma_{\gamma} = 1.4 $ as the integer $q$ in Eq. \eqref{dataProjection} is increased, with $L = 5$, $p = 5$, and $N = 2$.}
\label{q_var_high}
\end{table}

The effect of the accuracy of the data projection is also studied, increasing $q$ in Eq. \eqref{dataProjection} from $4$ to $7$ 
with $L=5$, $p=5$, and $N =2$. 
Results are listed in Tables \ref{q_var_low},
\ref{q_var_middle}, and \ref{q_var_high} for the small, medium, and large values of $\sigma_{\gamma}$, respectively (i.e. 0.08, 0.8, 1.4). 
We find that the number of the GMRES iterations for all methods is invariant of $q$,
suggesting that the rate of convergence is not influenced by an increasing density of the stiffness matrix. 
Similarly, the timing results are insensitive to the variation of $q$.  
The G-M-F-I and the G-M-Fg-I are comparable and outperform the other methods for the small and medium $\sigma_{\gamma}$, 
whereas the G-M-Fg-I performs the best for the large value of $\sigma_{\gamma}$.     
Note that the G-F-I is faster than the G-M-I for the small and medium $\sigma_{\gamma}$. 

\begin {table}[!ht] 
\begin{center}
	\begin{tabular}{|c|c|c|c|c|c|c|c|} \hline	
	 Precond. & {G-I} & {G-F-I} & {G-M-I} & {G-M-F-I} & {G-M-Fg-I} \\ \hline
	     $N$ & \shortstack{\\CPU Time [s] \\ (GMRES iter.)}  & \shortstack{\\CPU Time [s] \\ (GMRES iter.)}  & \shortstack{\\CPU Time [s] \\ (GMRES iter.)}  & \shortstack{CPU Time [s] \\ (GMRES iter.)} & \shortstack{CPU Time [s] \\ (GMRES iter.)} \\ \hline 	     
	     2 & \shortstack{\\4.74   \\ (129)} & \shortstack{\\2.82  \\(15)} & \shortstack{\\1.99  \\ (7)} & \shortstack{\\0.80  \\(6)} & \shortstack{\\0.71  \\(6)}\\ \hline 
	     3 & \shortstack{\\24.66  \\ (141)} & \shortstack{\\6.84  \\(15)} & \shortstack{\\12.12 \\ (7)} & \shortstack{\\3.43  \\(6)} & \shortstack{\\2.91  \\(6)}\\ \hline 
	     4 & \shortstack{\\102.61 \\ (145)} & \shortstack{\\15.77 \\(15)} & \shortstack{\\61.64 \\ (7)} & \shortstack{\\13.36 \\(6)} & \shortstack{\\11.57 \\(6)}\\ \hline 
	\end{tabular}
\end{center}
\caption{CPU times and GMRES iterations with $\sigma_{\gamma} = 0.08 $ as the stochastic input dimension $N$ is increased, with $L = 5$  and $p =q = 3$.}
\label{N_var_low}
\end{table}

\begin {table}[!ht] 
\begin{center}
	\begin{tabular}{|c|c|c|c|c|c|c|c|} \hline	
	 Precond. & {G-I} & {G-F-I} & {G-M-I} & {G-M-F-I} & {G-M-Fg-I}\\ \hline
	     $N$  & \shortstack{\\CPU Time [s] \\ (GMRES iter.)}  & \shortstack{\\CPU Time [s] \\ (GMRES iter.)}  & \shortstack{\\CPU Time [s] \\ (GMRES iter.)}  & \shortstack{CPU Time [s] \\ (GMRES iter.)} & \shortstack{CPU Time [s] \\ (GMRES iter.)}\\ \hline 	     
	     2 & \shortstack{\\6.53   \\ (190)} & \shortstack{\\3.93  \\(22)} & \shortstack{\\1.99  \\ (7)} & \shortstack{\\0.97  \\(8)} & \shortstack{\\0.79  \\(7)}\\ \hline 
	     3 & \shortstack{\\30.85  \\ (190)} & \shortstack{\\9.93  \\(23)} & \shortstack{\\12.12 \\ (7)} & \shortstack{\\4.19  \\(8)} & \shortstack{\\3.25  \\(7)}\\ \hline 
	     4 & \shortstack{\\139.39 \\ (230)} & \shortstack{\\24.72 \\(25)} & \shortstack{\\62.07 \\ (7)} & \shortstack{\\17.90 \\(9)} & \shortstack{\\12.96 \\(7)}\\ \hline 
	\end{tabular}
\end{center}
\caption{CPU times and GMRES iterations with $\sigma_{\gamma} = 0.8 $ as the stochastic input dimension $N$ is increased, with $L = 5$  and $p =q = 3$.}
\label{N_var_middle}
\end{table}

\begin {table}[!ht] 
\begin{center}
	\begin{tabular}{|c|c|c|c|c|c|c|c|} \hline	
	 Precond. & {G-I} & {G-F-I} & {G-M-I} & {G-M-F-I} & {G-M-Fg-I}\\ \hline
	     $N$ & \shortstack{\\CPU Time [s] \\ (GMRES iter.)}  & \shortstack{\\CPU Time [s] \\ (GMRES iter.)}  & \shortstack{\\CPU Time [s] \\ (GMRES iter.)}  & \shortstack{CPU Time [s] \\ (GMRES iter.)} & \shortstack{CPU Time [s] \\ (GMRES iter.)}\\ \hline 	     
	     2 & \shortstack{\\7.48   \\ (223)} & \shortstack{\\4.40  \\(25)} & \shortstack{\\1.98  \\ (7)} & \shortstack{\\1.32   \\(12)} & \shortstack{\\0.86   \\(8)}\\ \hline 
	     3 & \shortstack{\\34.06  \\ (215)} & \shortstack{\\12.25 \\(29)} & \shortstack{\\12.46 \\ (8)} & \shortstack{\\6.89   \\(15)} & \shortstack{\\3.57   \\(8)}\\ \hline 
	     4 & \shortstack{\\156.56 \\ (269)} & \shortstack{\\29.57 \\(30)} & \shortstack{\\61.92 \\ (8)} & \shortstack{\\34.82  \\(20)} & \shortstack{\\14.43  \\(8)}\\ \hline 
	\end{tabular}
\end{center}
\caption{CPU times and GMRES iterations with $\sigma_{\gamma} = 1.4 $ as the stochastic input dimension $N$ is increased, with $L = 5$  and $p =q = 3$.}
\label{N_var_high}
\end{table}

At last, we investigate the influence of the dimension of the stochastic input $N$, which influences the dimension of the stochastic finite dimensional space, given that $M_p = (N+p)!/(N! p!)$. 
We vary $N$ from $2$ to $4$ with $L=5$ and $p=q=3$. 
The numerical results are displayed in Tables \ref{N_var_low},
\ref{N_var_middle}, and \ref{N_var_high} for the small, medium, and large values of $\sigma_{\gamma}$, respectively, (i.e. 0.08, 0.8, 1.4).
For the small and medium $\sigma_{\gamma}$ in Tables \ref{N_var_low} and \ref{N_var_middle},  
the number of GMRES iterations remains stable when $N$ increases for G-F-I, G-M-I, G-M-F-I, and G-M-Fg-I, 
whereas it increases for G-I. 
Moreover, the performance of G-M-Fg-I is the best in terms of computational time.
For the large $\sigma_{\gamma}$ results are shown in Table \ref{N_var_high}, where the number of GMRES iterations increases for 
G-I, G-F-I, and G-M-F-I but it remains stable for both G-M-I and G-M-Fg-I ,
suggesting that the only rate of convergence that does not depend on $N$ is that of G-M-I and G-M-Fg-I.
Again we find that G-M-Fg-I is the most efficient in terms of CPU time. For instance if $N=4$ and  $\sigma_{\gamma}=1.4$, 
G-M-Fg-I is at least 2 times faster than G-M-F-I and G-F-I,
4 times faster than G-M-I, and 10 times faster than G-I.

For the stochastic diffusion problem with lognormal coefficient, after considering the effects of a variation of the mesh refinement, $\sigma_{\gamma}$, $p$, $q$ and $N$, 
we can hereby conclude that G-M-Fg-I is the best preconditioned solver among the ones considered in this work.  

\section{Conclusions}
A computational comparison of several preconditioners for the GMRES algorithm applied to the stochastic diffusion problem with lognormal coefficient has been presented.
The main focus of our work was to highlight the performances of field-split type preconditioners, that are inspired by physics based preconditioners for deterministic problems. A thorough analysis of the robustness and dependence on the mesh size and stochastic parameters showed that the field-split type preconditioners used within a multigrid framework are the best option, among the methods tested, especially if a specific grouping strategy is chosen.The best preconditioned solver presented (G-M-Fg-I) showed independence on mesh refinement, data projection (i.e. density of the matrix), input standard deviation, stochastic dimension (i.e. number of stochastic variables) and polynomial order of the SGM approximate solution. 
Moreover, it was always the fastest in terms of CPU time, which is of great importance when dealing with computationally heavy solutions of linear systems arising from a stochastic Galerkin type of discretization. These features make the G-M-Fg-I a highly valuable preconditioned solver for the specific problem at hand. The methodologies used in this work will be applied in the future to more complex stochastic PDEs to further investigate the potential of the field-split approach and, in particular, of the G-M-Fg-I preconditioned solver.
\bibliographystyle{plain}
\bibliography{UQ_Poisson}
\end{document}